\documentclass[preprint,12pt]{elsarticle}

\usepackage{amsthm,amsbsy,amsfonts,amssymb,amsmath,amscd,latexsym,mathabx}
\usepackage{fontenc}
\usepackage[mathscr]{euscript}
\usepackage{graphicx}
\usepackage{epic,eepicemu}
\usepackage{mathrsfs}
\usepackage{bm}
\usepackage[usenames]{color} 
\usepackage{xcolor} 
\definecolor{mygray}{gray}{0.85}
\definecolor{mygray-85}{gray}{0.85}
\definecolor{mygray-75}{gray}{0.70}
\definecolor{MyDarkBlue}{rgb}{0.2,0,1}
\usepackage[urlbordercolor = MyDarkBlue, urlcolor = MyDarkBlue]{hyperref}
 \usepackage[all]{xy} 
 \CompileMatrices
\newcommand{\markedrightarrow}[1]{\setlength{\unitlength}{1ex}
\begin{picture}(5,2)(0,0)
\put(2.5,1.5){\makebox(0,0)[c]{$\scriptstyle{#1}$}}
\put(2.5,0.4){\makebox(0,0)[c]{$\longrightarrow$}}
\end{picture}}

\newcommand{\myscripthat}[1]{
{\mbox{}}\;\widehat{\makebox[10pt][r]{$\mathscr{#1}$}}}

\newcommand{\mywidehat}[1]{
{\mbox{}}\;\widehat{\makebox[10pt][r]{${#1}$}}}

\newcommand{\dual}[1]{#1\raisebox{5.3pt}{$\scriptstyle{\vee}\!$}}
\newcommand{\duall}[1]{#1\hspace*{1pt}\raisebox{7pt}{$\scriptstyle{\vee}\!$}}

\newcommand{\proofsquare}{\setlength{\unitlength}{0.144ex}
\begin{picture}(10,10)(0,0)
\put(-0.26,0){\line(1,0){10.535}}
\put(-0.26,10){\line(1,0){10.535}}
\put(0,0){\line(0,1){10}}
\put(10,0){\line(0,1){10}}
\end{picture}}

\newtheorem{theorem}{Theorem}[section]
\newtheorem{proposition}[theorem]{Proposition}
\newtheorem{corollary}[theorem]{Corollary}
\newtheorem{lemma}[theorem]{Lemma}

\theoremstyle{definition}
\newtheorem{definition}[theorem]{Definition}

\theoremstyle{remark}
\newtheorem{remark}[theorem]{Remark}

\biboptions{square} 
\journal{Journal of Pure and Applied Algebra}

\begin{document}

\begin{frontmatter}



\title{A Closer Look at the Stacks of Stable Pointed Curves}


\author{Finn F. Knudsen}

\address{Department of Mathematical Sciences\\
Norwegian University of Science and Technology\\
Sentralbygg 2, Alfred Getz vei 1, 7491 Trondheim, Norway}
\begin{abstract}

In the theory of the moduli-stacks of stable $n$-pointed curves, there are two fundamental functors, {\emph{contraction}} and {\emph{stabilization}}. These functors are constructed in \cite{Knudsen:83}, where they are used to show that the various ${\overline{\mathcal{M}}}_{g,n}$'s are DM-stacks. In the treatment of stabilization,  there is a claim stated about the deformations of a pointed ordinary double point. This claim is not proved in \cite{Knudsen:83}. The goal of this paper is to prove this claim.  We decided not to use descent, so the equations are less simple than in \cite{Knudsen:83}, but the treatment is more elementary and the results seemingly more general. 
\end{abstract}


\end{frontmatter}

\section*{Introduction}
\label{intr}

Recall the definition of the stacks ${\overline{\mathcal{M}}}_{g,n}$ and ${\overline{\mathcal{C}}}_{g,n}$, for $n \geq 0$, $n+2g \geq 3$. They are both fibered categories over the category of all schemes. The objects are diagrams of the form
\begin{center}
\parbox[][9ex][c]{13ex}{$
\SelectTips{cm}{} 
\xymatrix{C \ar[d]_{{\pi}_{C}} \\ S \ar@/_3.5ex/[u]_{s_{1},\ldots,s_{n}}} $}
\quad \text{in ${\overline{\mathcal{M}}}_{g,n}$,  and of the form }\quad 
\parbox[][9ex][c]{13ex}{$
\SelectTips{cm}{} 
\xymatrix{C \ar[d]_{{\pi}_{C}} \\ S \ar@/_3.5ex/[u]_{s_{1},\ldots,s_{n}, \Delta}} $}\qquad \text{in ${\overline{\mathcal{C}}}_{g,n}$.}
\end{center}

The morphism ${\pi}_{C}$ is flat and proper. The morphisms $s_{1},\ldots,s_{n}$ are non-crossing sections not meeting any double points. There are no conditions on the extra section $\Delta$. The geometric fibers of ${\pi}_{C}$ are stable{\footnote{In  \cite{Knudsen:83} we use the term $n$-pointed stable curve, but as remarked by Frans Oort, it is not necessarily a stable curve with marked points.}} $n$-pointed curves, without and with one extra point respectively. An $n$-pointed curve is stable if it is reduced with at most ordinary double points, connected and has only a finite number of automorphisms fixing the $n$ distinct distinguished smooth points. Morphisms in these categories are cartesian diagrams commuting with the sections. The functor to the category of schemes takes an object to its base $S$. 

\smallskip

We show that ${\overline{\mathcal{C}}}_{g,n} $ is a universal stable $n$-pointed curve over ${\overline{\mathcal{M}}}_{g,n} $. For each object  $({\pi}_{C}; s_{1},\ldots,s_{n},\Delta)$ in ${\overline{\mathcal{C}}}_{g,n}$, we can simply forget the section $\Delta$, and we get an object  $({\pi}_{C}; s_{1},\ldots,s_{n})$ in ${\overline{\mathcal{M}}}_{g,n}$. This is a morphism of fibered categories, and we call it ${\sigma}_{n} :  {\overline{\mathcal{C}}}_{g,n} \rightarrow {\overline{\mathcal{M}}}_{g,n}$. For each object  $({\pi}_{C}; s_{1},\ldots,s_{n})$ in ${\overline{\mathcal{M}}}_{g,n}$, and for each $i$ with $1 \leq i \leq n$, we may define an extra section $\Delta = s_{i}$. This defines an object  $({\pi}_{C}; s_{1},\ldots,s_{n},\Delta)$ in ${\overline{\mathcal{C}}}_{g,n}$, so in this way we have $n$ sections  ${\sigma}_{n,i} : {\overline{\mathcal{M}}}_{g,n} \rightarrow {\overline{\mathcal{C}}}_{g,n}$. Note that if ${\pi}_{C} : C \rightarrow S$ with sections $s_{i} : S \rightarrow C$, $1 \leq i \leq n$ is an object of ${\overline{\mathcal{M}}}_{g,n}$, the map ${\mathrm{pr}}_{2} : C\times_{S} C \rightarrow C$ with sections $s_{n,i} = (s_{i} \circ {\pi}_{C}, \bm{1}_{C})$ and the extra section $\Delta = (\bm{1}_{C}, \bm{1}_{C})$ is an object of ${\overline{\mathcal{C}}}_{g,n}$. Identifying the schemes $S$ and $C$ with their respective slice categories ${\mathrm{Schemes}}/S$ and ${\mathrm{Schemes}}/C$, the reader should convince himself or herself that we have a cartesian commutative diagram of fibered categories

\smallskip

\begin{center}
\parbox[c][9ex][c]{27ex}{$
\SelectTips{cm}{} 
\xymatrix{C \ar[d]_{{\pi}_{C}}   \ar[r]   & {\overline{\mathcal{C}}}_{g,n}    \ar[d]_{{\sigma}_{n}}    \\  S  \ar@/^3.5ex/[u]^{s_{1},\ldots,s_{n}}  \ar[r]  & {\overline{\mathcal{M}}}_{g,n} \ar@/_3.5ex/[u]_{\sigma_{n,1},\ldots, \sigma_{n,n}} \, ,} $}
\end{center}

\smallskip

which shows that ${\sigma}_{n}$ is representable and that $({\sigma}_{n}; {\sigma}_{n,1}, \ldots , {\sigma}_{n,n})$ is a universal stable $n$-pointed curve.

\smallskip

As mentioned in the abstract, there is a morphism $c : {\overline{\mathcal{M}}}_{g,n{+}1} \rightarrow {\overline{\mathcal{C}}}_{g,n}$ called {\emph{contraction}} and a morphism $s : {\overline{\mathcal{C}}}_{g,n} \rightarrow {\overline{\mathcal{M}}}_{g,n{+}1}$ called {\emph{stabilization}}. These morphisms are constructed in \cite{Knudsen:83}, where it is proved that they are inverse to each other. Since ${\sigma}_{n}$ is representable, so is ${\pi}_{n{+}1} = {\sigma}_{n} c$ and consequently  ${\pi}_{n{+}1} :  {\overline{\mathcal{M}}}_{g,n{+}1} \rightarrow {\overline{\mathcal{M}}}_{g,n}$ is also a universal stable $n$-pointed curve. This is then used to show inductively in $n$ that the various $ {\overline{\mathcal{M}}}_{g,n} $'s are DM-stacks.

\smallskip

The construction of the stabilization morphism and the proof that contraction and stabilization are inverse to each other is a bit tricky. The way stabilization is constructed is as follows. For every object  $({\pi}_{C}, s_{1},\ldots,s_{n},\Delta)$ in ${\overline{\mathcal{C}}}_{g,n}$, where ${\pi}_{C} : C \rightarrow S$, we define the coherent sheaf ${\mathscr{K}}_{({\pi}_{C}, s_{1},\ldots,s_{n},\Delta)} = {\mathscr{K}}$ as the cokernel of the map $\delta$ in the short exact sequence
\begin{align}\label{fund:exactsequence}
&0 \rightarrow  {\mathscr{O}}_{C}  \markedrightarrow{\delta} \dual{\mathscr{J}} \oplus {\mathcal{O}}_{C} (s_{1}+s_{2}+ \cdots +s_{n})  \markedrightarrow{p}  {\mathscr{K}} \rightarrow 0,
\end{align}
where ${\mathscr{J}}$ is the defining ideal of the section $\Delta$ and $\delta$ is the ``diagonal" $\delta(t) = (t,t)$. We define $C^{s} = {\textbf{Proj}}_{C}({\mathrm{Sym}}\hspace*{1pt}{\mathscr{K}})$. The union of the images of the sections $s_{i}$ and $\Delta$ support sheaves ${\mathscr{L}}_{s} = {\mathscr{K}}/p(\dual{\mathscr{J}}\; )$ and ${\mathscr{L}}_{\Delta} = {\mathscr{K}}/p({\mathcal{O}}_{C} (s_{1}+s_{2}+ \cdots +s_{n}))$, and there are surjections $s_{i}^{*}{\mathscr{K}} \rightarrow s_{i}^{*}{\mathscr{L}}_{s}$ and ${\Delta}^{*}{\mathscr{K}} \rightarrow {\Delta}^{*}{\mathscr{L}}_{\Delta}$. The important fact is that (\ref{fund:exactsequence}) commutes with base-change and that the sheaves $s_{i}^{*}{\mathscr{L}}_{s}$ and ${\Delta}^{*}{\mathscr{L}}_{\Delta}$ are invertible. Hence these surjections define liftings of the sections making $C^{s}$ together with the lifted sections a stable $(n{+}1)$-pointed curve. We call it the stabilized curve.

\smallskip

The hard part is to prove these facts in the case where the section ${\Delta}$ passes through a double point of a fiber. This is the only case where ${\mathscr{J}}$ is not an invertible sheaf; nevertheless it has certain good properties, stated in \cite{Knudsen:83} Section 2, Lemma 2.2. It is our Main Lemma stated below. To prove this lemma we apply the deformation theory of pointed double points, Proposition \ref{prop:versal}, to get explicit descriptions of formal neighborhoods of $k$-valued points in ${\overline{\mathcal{C}}}_{g,n}$. The Main Lemma then follows from the Key Example using the properties of the concept of {\emph{stable reflexivity}}{\footnote{This is the property that ensures the compatibility with base change that was sought earlier.}}, which is proved in  \cite{Knudsen:83}, Appendix. Recently Runar Ile  \cite{Ile:11} has given another characterization of this concept, together with another proof of a very nice generalization of the Main Lemma. In the final section we complete the family of stabilized curves along the fiber curve and compare it to the formal scheme coming from the Key Example in order to prove stability.

\medskip

$ \mbox{}$ \hspace*{-4.5ex} {\bf{Main Lemma.}} {\emph{
Let $S$ be a noetherian scheme, and $\pi : C \rightarrow S$ a flat family of curves whose geometric fibers are reduced curves with at most ordinary double points. Let $\Delta : S \rightarrow C$ be a section defined by an ${\mathcal{O}}_{C}$-ideal ${\mathscr{J}}$. Then 
\begin{itemize}
\item[$1)$] The sheaf  ${\mathscr{J}}$ is stably reflexive with respect to $\pi$.
\item[$2)$] Let ${\mathscr{J}}'$ be the subsheaf of the total quotient ring sheaf $K_{C}$ consisting of sections that multiply ${\mathscr{J}}$ into ${\mathcal{O}}_{C}$ and let $\dual{\mathscr{J}} = \;${\emph{\textbf{Hom}}}$_{{\mathcal{ \hspace*{1pt} O}}_{C}}({\mathscr{J}},{\mathcal{O}}_{C})$. Then ${\mathscr{J}}' \approx \dual{\mathscr{J}}$.\item[$3)$] The pullback ${\Delta}^{\ast}(\dual{\mathscr{J}}/{\mathcal{O}}_{C})$ is an invertible sheaf on $S$.
\end{itemize}}}

\section{Pointed ordinary double points}
We consider curves over a base consisting of a single reduced point ${\mathrm{Spec}}(k)$, where  $k$ is a field{\footnote{In \cite{Knudsen:83} we work in the {\'{e}}tale topology. The Main Lemma then follow by descent. In this note however we will stick to the Zariski topology, so even at an ordinary double point
that is itself rational, its tangent lines do not have to be so. It turns out that we get slightly less simple equations to work with. The reader may, if he or she wishes to, consider only a case where we have rational tangents, e.g. $\gamma = 1$, $\delta = 0$. My feeling is that the proofs are almost identical, so not very much is gained by doing so.}}. For any $n \in {\mathbb{N}}$ and any element $f \in k[[x,y]]$ let $f_{n}$ be the $n$ th homogeneous part of $f$. We denote by  $S_{n}$ the vector space of homogeneous elements of degree $n$  in $k[x,y] \subseteq  k[[x,y]]$.
\begin{definition}\label{def:approx} \quad With notation as above, let $q(x,y) = x^{2} + \gamma xy + \delta y^{2} \in S_{2}$.  For any $n \geq 0$, the map $Q_{n}: S_{n} \times S_{n} \rightarrow S_{n{+}1}$ is the $(n{+}1)$st.{\hspace*{-0.2ex}} homogeneous part of the polynomial $q(x+\mu , y+\nu) - q(x,y)$. In other words, $Q_{n}(\mu,\nu) =  (2 \mu + \gamma \nu)x + (\gamma \mu + 2 \delta \nu)y$. 
\end{definition}
\begin{lemma}\label{lem:approx} \quad Let $q \in S_{2}$ be a quadratic form. The associated map $Q_{n}: S_{n} \times S_{n} \rightarrow S_{n{+}1}$ is linear, and $q(x+\mu , y+\nu) - q(x,y) - Q_{n}(\mu, \nu) \in S_{2n}$. If $q$ is non-degenerate, $Q_{n}$ is surjective for all $n \geq 0$.
\end{lemma}
\begin{proof}\quad We have $Q_{n}(\mu, \nu) = (2 \mu + \gamma \nu)x + (\gamma \mu + 2 \delta \nu)y$. Any $f_{n{+}1} \in S_{n{+}1}$ can be written $f_{n{+}1} = u_{n}x + v_{n}y$. If we define 
${\mu}_{n} = -2\delta  u_{n} + \gamma v_{n}$ and ${\nu}_{n} = \gamma  u_{n} - 2 \delta v_{n}$, a calculation shows that $Q_{n}({\mu}_{n}, {\nu}_{n}) = ({\gamma}^{2} - 4 \delta) f_{n{+}1}$. The lemma follows since that $q$ is non-degenerate means that ${\gamma}^{2} - 4 \delta$ is invertible.
\end{proof}
\begin{proposition}\label{prop:odp} \quad Let $k$ be a field, $\pi: C \rightarrow {\mathrm{Spec}}(k)$ and $\Delta : {\mathrm{Spec}}(k) \rightarrow C$ be such that the geometric fibers of $\pi$ are reduced curves with at most ordinary double points and assume that  $c=\Delta({\mathrm{Spec}}(k))$ is an ordinary double point. Let $(R,{\bm{m}})$ be the local ring of $c \in C$. Then there are elements $\gamma$ and $\delta$ in $k$ with ${\gamma}^{2} - 4\delta$ not equal to zero such that $\widehat{R} \approx  k[[x,y]]/(x^{2} + \gamma xy + \delta y^{2})$.
\end{proposition}
\begin{proof}\quad Since ${\bm{m}}/{\bm{m}}^{2} \otimes_{k} {\overline{k}}$ is a vector space of dimension 2, so is ${\bm{m}}/{\bm{m}}^{2}$. If we choose two elements of ${\bm{m}}$ that generate ${\bm{m}}/{\bm{m}}^{2}$ we get a surjection $k[[x,y]] \rightarrow \widehat{R}$. The morphism $\pi$ is a locally complete intersection morphism, so the kernel is generated by a single power series $f=\sum_{n=2}^{\infty} f_{n}$. After possibly a making a linear change of variables and multiplying $f$ by a suitable scalar, we may assume that $f_{2}(x,y) = q(x,y) = x^{2} + \gamma xy + \delta y^{2}$. Since at the unique point in the geometric fiber over the point $c$ the curve has two branches with distinct tangents, the quadratic form $q$ is non-degenerate. Let $Q_{n}$, be as in Lemma \ref{lem:approx}, and for each $n \geq 0$ let $P_{n{+}1}: S_{n{+}1} \rightarrow S_{n} \times S_{n}$ be a right inverse.

\smallskip

For all $n \geq 1$, we define homogeneous polynomials  ${\epsilon}_{n{+}1}$, ${\mu}_{n}$, ${\nu}_{n}$, $x_{n}(x,y)$ and $y_{n}(x,y)$  of degree $n$, by recursion. For $n=1$ put ${\epsilon}_{2}=0$, ${\mu}_{1}=0$, ${\nu}_{1}=0$, $x_{1}(x,y) = x$ and $y_{1}(x,y)=y$. Having defined ${\epsilon}_{n{+}1}$, ${\mu}_{n}$, ${\nu}_{n}$, $x_{n}(x,y)$ and $y_{n}(x,y)$, we define  
\begin{align*}
&{\epsilon}_{n+2} = \text{ the $(n{+}2)$nd.{\hspace*{-0.2ex}} homogeneous component of $q(x_{n}, y_{n}) - f(x,y)$, } \\
&({\mu}_{n{+}1}, {\nu}_{n{+}1}) = -P_{n+2}({\epsilon}_{n+2}) ,\\
& x_{n{+}1} = x_{n} + {\mu}_{n{+}1}, \\
& y_{n{+}1} = y_{n} + {\nu}_{n{+}1}. 
\end{align*}

\smallskip

We claim that for all $n \geq 1$, the order of $q(x_{n}, y_{n}) - f(x,y)$ is at least $n+2$. The claim certainly holds for $n=1$. Assume that $q(x_{n}, y_{n}) - f(x,y) = \epsilon_{n+2} + h(x,y)$, with the order of $h$ at least $n+3$. We have 
\begin{align*}
&q(x_{n{+}1}, y_{n{+}1}) - f(x,y) = \\
&q(x_{n} + \mu_{n{+}1}, y_{n} + \nu_{n{+}1}) - q(x_{n}, y_{n}) + q(x_{n}, y_{n}) - f(x,y) = \\
&Q_{n{+}1}(\mu_{n{+}1}, \nu_{n{+}1}) + \epsilon_{n+2} + h(x,y) + q( \mu_{n{+}1},  \nu_{n{+}1}) =  \\
&h(x,y) + q( \mu_{n{+}1},  \nu_{n{+}1}).
\end{align*}The sequences $\{x_{n}\}_{n \geq 1}$ and $\{y_{n}\}_{n \geq 1}$ are Cauchy sequences, and if we define $x'= \lim_{n \rightarrow \infty} x_{n}$ and $y'= \lim_{n \rightarrow \infty}y_{n}$, $\,(x,y) \mapsto (x', y')$ is a continuous change of coordinates, and $\widehat{R}= k[[x',y']]/(x'^{2} + \gamma x'y' + \delta y'^{2})$.
\end{proof}
\section{The deformation theory of a pointed ordinary double point} 
In this section we prove, in a slightly more general form, the claim in the middle of the proof of Lemma 2.2 on page 175 in \cite{Knudsen:83}. This claim is a special case of Proposition \ref{prop:versal} on the next page. All concepts in this section can be found in \cite{Schlessinger:68}. Let $\Lambda$ be a noetherian complete local ring with maximal ideal ${\bm{\mu}}={\bm{\mu}}_{\Lambda}$ and quotient field $k=\Lambda / {\bm{\mu}}$.  We define ${\bm{C}}={\bm{C}}_{\Lambda}$ to be the category of maps $A \rightarrow k$, where $A$ is an Artinian local $\Lambda$-algebras with residue field $k$. The category $\widehat{\bm{C}}=\widehat{\bm{C}}_{\Lambda}$ will denote the category of maps $A \rightarrow k$, where $A$ is a complete noetherian local $\Lambda$-algebra for which $A/{\bm{m}}^{n} \rightarrow k$ is in  ${\bm{C}}$ for all $n$.

\smallskip

Let  $\gamma , \delta \in k$ be elements such that ${\gamma}^{2} - 4\delta \neq 0$, and let  $q(X,Y) = X^{2} + \gamma XY + \delta Y^{2}$. We define $R_{0}=k[[X,Y]]/(q(X,Y))$ and denote by  $u$ and $v$ the classes of $X$ and $Y$ respectively.

\smallskip

We consider the category ${\mathscr{F}}$ with objects commutative cocartesian diagrams of local rings of the form shown in the figure below, and such that $R$ is a flat $A$ algebra, complete with respect to the linear topology defined by the maximal ideal ${\bm{m}}_{R}$. All maps are local homomorphisms.
\begin{center}
$
\SelectTips{cm}{} 
\xymatrix{R \ar[r]  \ar@/_3ex /[d]_{{\Delta}} & R_{0}  \ar@/^3ex /[d]^{{\Delta}_{0}}\\ A \ar[r]  \ar[u] & k \ar[u]  } 
$
\end{center}
We define the category  $\myscripthat{F}$ over $\widehat{\bm{C}}$ similarly. Note that ${\mathscr{F}}$ and ${\bm{C}}$ are full subcategories of 
$\myscripthat{F}$ and $\widehat{\bm{C}}$ respectively. The category ${\mathscr{F}}$ is cofibered in groupoids over  ${\bm{C}}$. For each object $A \rightarrow k$ the fiber category ${\mathscr{F}}_{A}$ is a groupoid. We let $F(A)$ denote the set of isomorphism-classes in ${\mathscr{F}}_{A}$, and by  $\widehat{F}(A)$ the set of isomorphism-classes in  $\;\widehat{\makebox[10pt][c]{$\mathscr{F}_{A}$}}$. Of course $F$ is the restriction of $\widehat{F}$ to ${\bm{C}}$.
\begin{proposition} \label{prop:versal}
Let ${\tilde{\gamma}}$ and  ${\tilde{\delta}}$ be any liftings of $\gamma$ and $\delta$ to $\Lambda$, and let ${\tilde{q}}(X,Y)= X^{2} + {\tilde{\gamma}} XY + {\tilde{\delta}} Y^{2}$. Define the complete noetherian local rings $A_{pd} = {\Lambda}[[S,T]]$, and $R_{pd} = A_{pd}[[X,Y]]/({\tilde{q}}(X,Y)-{\tilde{q}}(S,T))$. We denote by $U$ and $V$ the classes of $X$ and $Y$ in $R_{pd}$ respectively. Then $R_{pd}$ is a flat $A_{pd}$-algebra. We have a local homomorphism ${\Delta}_{pd} : R_{pd} \rightarrow A_{pd}$ and $f: R_{pd} \rightarrow R_{0}$ defined by ${\Delta}_{pd}(X)=S$, ${\Delta}_{pd} (Y)=T$, $f(U)=u$, $f(V)=v$ and $f(S) = f(T) = 0$. The diagram 
\begin{align}
\xi &= \parbox[][3ex][r]{23ex}{$\SelectTips{cm}{} 
\xymatrix{R_{pd} \ar[r]^f  \ar@/_3ex /[d]_{{\Delta}_{pd} } & R_{0}  \ar@/^3ex /[d]^{{\Delta}_{0}}\\ A_{pd} \ar[r]  \ar[u] & k \ar[u] }$}
\end{align}
is an object of  \, $\widehat{\makebox[7pt][r]{$\mathscr{F}$}}$, and up to isomorphism it does not depend on the liftings. Moreover this object is a hull for the functor $F$.
\end{proposition}
We prove the proposition by first studying the coslice categories $\myscripthat{G} = (A_{pd}\, {\downarrow} \,\widehat{\bm{C}})$ and ${\myscripthat{H}} = (\xi \, {\downarrow }\myscripthat{F})$. We denote the restriction of these categories to $\bm{C}$, by $\mathscr{G}$ and $\mathscr{H}$ respectively. By abstract nonsense $\myscripthat{H}$ is cofibered in groupoids over $\myscripthat{G}$, and both categories are cofibered in groupoids over $\widehat{\bm{C}}$. We denote by $\widehat{G}(A)$ and  $\widehat{H}(A)$ the set of isomorphism classes in the fiber of an object $A \rightarrow k$ in $\widehat{\bm{C}}$. Both $\widehat{H}$ and $\widehat{G}$ are covariant functors to sets and the forgetful functor induces an isomorphism $\widehat{H} \rightarrow \widehat{G}$ which restricts to an isomorphism between their restrictions to $\bm{C}$, $H \rightarrow G$. By definition these functors are pro representable. Then we prove that the forgetful functor $\mathscr{H}\rightarrow \mathscr{F}$ is smooth \cite{Schlessinger:68} Definition 2.2, and induces an isomorphism between the tangent spaces. This will prove the proposition.

\begin{definition}\label{def:quadraticform} 
For any object $A \rightarrow k$ in $\widehat{\bm{C}}$ we denote by ${\gamma}_{A}$ and ${\delta}_{A}$, the images of the elements ${\tilde{\gamma}}$ and ${\tilde{\delta}}$ by the structure map. We denote by $q_{A}$ the quadratic form $q_{A}(X,Y) = X^{2} + {\gamma}_{A} XY + {\delta}_{A}Y^{2}$.
\end{definition}
\begin{definition}\label{def:planar} 
We will call any image of the object $\xi$ defined in Proposition \ref{prop:versal} {\emph{planar}}. A planar object occurs as the right square in a diagram, where the whole square and the left square are cocartesian.
\begin{center}
$
\SelectTips{cm}{} 
\xymatrix{R_{pd} \ar@/^3.4ex /[rr]^f  \ar@/_2ex /[d]_{{\Delta}_{pd} } \ar[r]^{g_{R}} &R  \ar@/^2ex /[d]^{{\Delta}_{R}} \ar[r]^{h_{R}} &R_{0}  \ar@/^2ex /[d]^{{\Delta}_{0}}\\ A_{pd} \ar@/_3.4ex /[rr]  \ar[u] \ar[r]^{g_{A}} &  A \ar[u] \ar[r]^{h_{A}}& k \ar[u] }
$

$\mbox{}$

$\mbox{}$

\end{center}
By abstract nonsense, the right square is also cocartesian, and hence an object of $\myscripthat{F}$, the whole diagram is an object of $\myscripthat{H}$, and the bottom line is object of $\myscripthat{G}$.
\end{definition}

\begin{remark} If an object 
\begin{center}
$
\SelectTips{cm}{} 
\xymatrix{R \ar[r]  \ar@/_3ex /[d]_{{\Delta}} & R_{0}  \ar@/^3ex /[d]^{{\Delta}_{0}}\\ A \ar[r]  \ar[u] & k \ar[u] }
$
\end{center}
of $\myscripthat{F}$ is planar, there are elements $u_{R}$ and $v_{R}$ mapping to $u_{R_{0}}$ and $v_{R_{0}}$ respectively, elements $s_{A}$ and $t_{A}$ in ${\bm{m}}_{A}$, such that  $q_{A}(u_{R}, v_{R}) = q_{A}(s_{A}, t_{A})$. We do not change names of elements in $A$ and their images in $R$.
\end{remark}

\begin{lemma}\label{lem:generatedunique}
Let $(A, {\bm{m}}_{A})$ be a complete local ring, and let $\gamma, \delta \in A$, $s,t \in  {\bm{m}}_{A}$ be arbitrary elements and let $q(X,Y)= X^{2} + \gamma XY + \delta Y^{2}$ be a quadratic form. Let $I$ be the ideal in $A[[X,Y]]$ generated by $q(X,Y) - q(s,t)$. Let $R = A[[X,Y]]/I$, and let $u$ and $v$ be the classes of $X$ and $Y$ in $R$ respectively. Then $R$ is flat over $A$ and every element of $R$ has a representation of the form $r=f(v) + ug(v)$, for unique power series  $f,g \in A[[Y]]$.
\end{lemma}
\begin{proof}
We define polynomials $f_{n}$, $g_{n}$ in $A[[Y]]$, $h_{n}$ in $A[[X,Y]]$ by recursion as follows. 
\begin{align*}
& f_{0} = g_{0} = h_{0} = 1, \\
& f_{1} = 0,\; g_{1} = -\gamma Y, \;h_{1} =  X + \gamma Y ,\\
& f_{2} = q(s,t) - \delta Y^{2},
\end{align*}
and for $n \geq 2$ 
\begin{align*}
f_{n{+}1} &= f_{2}g_{n{-}1}, \\
g_{n{+}1} &= f_{n{+}1} + g_{1}g_{n}, \\
h_{n{+}1} &= g_{n{+}1} + Xh_{n}. \\
\end{align*}
These formulas ensure that $f_{n}$, $g_{n}$ and $h_{n}$ all belong to the $n$th. power of the defining ideal ${\bm{m}}_{A[[X,Y]]}$, and induction shows that,
$$
 X^{n} = f_{n} + X g_{n{-}1} + h_{n{-}2} (q(X,Y) - q(s,t)).
$$

If $r = \sum a_{n} X^{n}$ is any power series in $A[[Y]][[X]]$, the series $f = \sum a_{n } f_{n}$, $g = \sum a_{n{+}1 } g_{n}$ and $h = \sum a_{n+2 } h_{n}$ all converge and
$$
 r = f + X g + h(q(X,Y) - q(s,t)).
$$
If $0 = f + X g + h(q(X,Y) - q(s,t))$, with $f,g$ in $A[[Y]]$, comparing $X$-degrees shows that $h = 0$ and therefore both $f=0$ and $g=0$. It follows that as an $A$-module $R \approx A[[Y]] \oplus A[[Y]]$ which is flat.
\end{proof}
\begin{lemma}\label{lem:newcoordinates}
Let $A$ be any commutative ring, and let $\tau$ be an element of $A$ with ${\tau}^{2} =0$. Let $q(X,Y) = X^{2} + \gamma XY + \delta  Y^{2}$ be a quadratic form with discriminant $d= {\gamma}^{2} - 4 \delta$ a unit in $A$, and let $f$ be a power series with vanishing constant term. Then there exist power series $\mu$ and $\nu$ in $A[[X,Y]]$ such that if we define $X' = X+{\tau}\mu$ and $Y'= Y+ {\tau}\nu$, then 
\begin{align}\label{eq:newcoordinates}
X^{2} + \gamma XY + \delta Y^{2} + {\tau}f(X,Y) &=   X'^{2} + \gamma X'Y' + \delta Y'^{2}.
\end{align}
Note that $X' = X+{\tau}\mu$ and $Y'= Y+ {\tau}\nu$ is a continuous change of variables because $X = X' - {\tau }g(X',Y')$ and $Y= Y'- {\tau }h(X',Y')$.
\end{lemma}
\begin{proof}\quad We may write $f= d(Xu + Yv)$. If we define $\mu = -2\delta u + \gamma v$ and $\nu = \gamma u - 2v$, a calculation shows that $X^{2} + \gamma XY + \delta Y^{2} + {\tau}f(X,Y) =   X'^{2} + \gamma X'Y' + \delta Y'^{2}$.
\end{proof}
\begin{lemma}\label{lem:allaregenerated}
Consider an object of the category  $\mathscr{F}$, as given by the diagram below.
\begin{center}
$
\SelectTips{cm}{} 
\xymatrix{R \ar[r]  \ar@/_3ex /[d]_{{\Delta}} & R_{0}  \ar@/^3ex /[d]^{{\Delta}_{0}}\\ A \ar[r]  \ar[u] & k \ar[u] }
$
\end{center}
Then if ${\tilde{u}}$ and ${\tilde{v}}$ are liftings of $u$ and $v$ respectively, every element $r \in R$ can be expressed as $r = f({\tilde{v}}) + {\tilde{u}}g({\tilde{v}})$, where and $f,g \in A[[Y]]$. 
\end{lemma}
\begin{proof}
Let $R' \subseteq R$ be the $A$-submodule generated by elements of the form $r = h({\tilde{v}}) + {\tilde{u}}k({\tilde{v}})$. Tensoring with $k$ we have a right-exact sequence
$$
R' \otimes_{A} k \rightarrow R \otimes_{A} k \rightarrow R/R' \otimes_{A} k \rightarrow 0.
$$
By Lemma \ref{lem:generatedunique} the first arrow is a surjection, so $R = R' + {\bm{m}}_{A}R$. Substituting this expression for $R$ back and taking into account that ${\bm{m}}_{A}R' \subseteq R'$, it follows that $R = R' + {\bm{m}}_{A}^{2}R$. By induction it follows that $R = R' + {\bm{m}}_{A}^{n}R$ for all $n$ and the lemma follows because there is an $n$ such that ${\bm{m}}_{A}^{n}=0$.
\end{proof}
The next lemma shows that ${\mathscr{H}} \rightarrow {\mathscr{F}}$ is smooth and hence essentially surjective.
\begin{lemma}\label{lem:smallisplanar}
Let $g_{A} : B \rightarrow A$ be a small map in ${\bm{C}}$ with ${\mathrm{ker}}\,g_{A}  = (\tau)$. Consider the diagram below which represents a morphism in  ${\mathcal{F}}$. 

\vspace*{-2ex}

\begin{center}
$
\SelectTips{cm}{} 
\xymatrix{S \ar@/^3.4ex /[rr]^{h_{S}} \ar@/_2ex /[d]_{{\Delta}_{S} } \ar[r]^{g_{R}} &R  \ar@/^2ex /[d]^{{\Delta}_{R}} \ar[r]^{h_{R}} &R_{0}  \ar@/^2ex /[d]^{{\Delta}_{0}}\\ B \ar@/_3.4ex /[rr]  \ar[u] \ar[r]^{g_{A}} &  A \ar[u] \ar[r]^{h_{A}}& k \ar[u] }
$
\end{center}
\bigskip
If the right square is planar, then so is the big square.
\end{lemma}
\begin{proof} 
Pick ${u_{S}}$ and $v_{S}$ in $S$ such that $g_{R}(u_{S})=u_{R}$ and $g_{R}(v_{S})=v_{R}$ and let ${\Delta}_{S}(u_{S})=s_{B}$ and ${\Delta}_{S}(v_{S})=t_{B}$. Then $g_{A}(s_{B})=s_{A}$ and $g_{A}(t_{B})=t_{A}$. Note also that if we alter ${u_{S}}$ and ${v_{S}}$ to ${u_{S}}'$ and ${v_{S}}'$ by elements of the form ${\tau}(a{u_{S}}+b{v_{S}})$, then ${s_{B}}= {\Delta}_{S}({u_{S}}')$ and  ${t_{B}}= {\Delta}_{S}({v_{S}}')$ because ${\tau}{\bm{m}}_{B} = \{0\}$. We have $g_{R}(q_{S}(u_{S}, v_{S}) - q_{S}(s_{B}, t_{B})) =0$, so by Lemma \ref{lem:allaregenerated} there is a power series $f$ with coefficients in $B$ such that $q_{S}(u_{S}, v_{S}) - q_{S}(s_{B}, t_{B}) = \tau f(u_{S}, v_{S})$. Since ${\Delta}_{S}(q_{S}(u_{S}, v_{S}) - q_{S}(s_{B}, t_{B})) = 0$, we have $0 = {\Delta}_{S}(\tau f) = \tau {\Delta}_{S}(f_{0}) = \tau f_{0}$. Therefore we have shown that $f$ is without constant term. By Lemma \ref{lem:newcoordinates} we can alter ${u_{S}} $ and ${v_{S}}$ such that $q_{S}(u_{S}, v_{S}) - q_{S}(s_{B}, t_{B}) = 0$. Note that the corresponding alterations $s_{B}' = s_{B} + \tau \Delta (\alpha)$ and $t_{B}' = t_{B} + \tau \Delta (\beta)$ do not alter $q_{S}(s_{B}, t_{B})$, since $q_{S}(s_{B}', t_{B}') - q_{S}(s_{B}, t_{B}) = \tau  \Delta (\alpha) \frac{\partial{q}}{\partial{x}}(s_{B}, t_{B}) +  \tau  \Delta (\beta) \frac{\partial{q}}{\partial{y}}(s_{B}, t_{B})$, and both partial derivatives are linear forms in $s_{B}$ and  $t_{B}$, hence are contained in the maximal ideal of $B$, and so the difference vanishes.
We need to prove that there are no more relations, so let $I$ be the kernel of the map $B[[x,y]]/(q_{B}(x,y) - q_{B}(s_{B}, t_{B})) \rightarrow S$. Tensoring with $A$ we get a long exact sequence
$$
{\mathrm{Tor}}^{B}(S, A) \rightarrow I \otimes_{B}A \rightarrow B[[x,y]]/(q_{B}(x,y) - q_{B}(s_{B}, t_{B}))\otimes_{B}A  \rightarrow S\otimes_{B}A \rightarrow 0
$$
Now $S$ is flat over $B$ and the last map is an isomorphism therefore $ I \otimes_{B}A \approx I/ \tau I= 0$ and since ${\tau}^{2}=0$ we have $I=0$.
\end{proof}
\begin{lemma} \label{lem:tangent}
The tangent map ${\bm{t}}_{H} \rightarrow {\bm{t}}_{F}$ is an isomorphism, and both spaces are isomorphic to the same vector space $k^{2}$.
\end{lemma}
\begin{proof}
Let $k[\epsilon]$ be the ring of dual numbers. Then since ${\epsilon}^{2} = 0$, every object of ${\mathscr{F}}(k[\epsilon])$ is of the form 
$$\SelectTips{cm}{} 
\xymatrix{R \ar[r]  \ar@/_3ex /[d]_{{\Delta}_{R}} & R_{0}  \ar@/^3ex /[d]^{{\Delta}_{0}}\\
k[\epsilon]\ar[r]  \ar[u] & k \ar[u] } 
$$
with $R \approx k[\epsilon, x, y]/(q(x, y), {\epsilon}^{2})$. Since ${\epsilon}^{2} = 0$, there is a well defined pair  $(a, b) \in k^{2}$ such that the isomorphism class of this object in ${\bm{t}}_{F}$ is determined by the pair $({\Delta}_{R}(x), {\Delta}_{R}(y)) = (\epsilon a, \epsilon b)$.

The space ${\bm{t}}_{H} = {\mathrm{Hom}}_{\bm{C}}({\Lambda}[[S, T]], k[\epsilon])$. Elements of this space are homomorphisms determined by the values of the elements $S$ and $T$, which is also a pair $(\epsilon a, \epsilon b)$ with $(a, b) \in k^{2}$. With these identifications, the map ${\bm{t}}_{H} \rightarrow {\bm{t}}_{F}$ induces the identity on $k^{2}$. This proves the lemma and Proposition \ref{prop:versal}.
\end{proof}

%\newpage

\section{The Key Example}\label{sec:key}
Let $A$ be a unitary noetherian commutative ring, $B= A[X,Y]$, $q(X,Y) = X^{2} + \gamma XY + \delta Y^{2}$ a quadratic form such that the discriminant ${\gamma}^{2} - 4\delta $ is a unit in $A$. Let $s$ and $t$ be arbitrary elements of $A$. The element $x =  q(X,Y) - q(s,t)$ is not a zero-divisor in $B$ and the reader may check that we have a matrix-factorization $ \phi \psi = \psi\phi = x{\cdot} {\bm{1}}_{B^{2}} $, where
$$ \phi=\begin{pmatrix}  \delta Y + \delta t + \gamma X & X+s + \gamma t \\ -(X-s )  &  Y-t \end{pmatrix} \quad \text{and}\quad   \psi=\begin{pmatrix} Y-t&-(X+s + \gamma t) \\X-s & \delta Y + \delta t + \gamma X \end{pmatrix}.$$
This matrix-factorization comes about from trying to find a minimal resolution of the ideal generated by $X-s$ and $Y-t$ modulo $q(X,Y) - q(s,t)$. If we define \quad $p= \begin{pmatrix} 0&-1\\1&0 \end{pmatrix}$, \quad we have\quad 
$ p \psi = {\phi}^{\text{tr}} p  \quad \text{and} \quad  p \phi = {\psi}^{\text{tr}} p$.
We let $R=B/(x)$ be the quotient ring and denote by $u$, $v$, $\alpha$, $\beta$ the classes of $X$, $Y$, $\phi$,  $\psi$ respectively. Then if we put $E = R \oplus R$ it follows from  \cite{Eisenbud:80} Proposition 5.1 page 49 that the diagram below commutes and that the rows are exact.
$$\SelectTips{cm}{} 
\xymatrix@1{\ar[r]^{\alpha} &\:E \ar[r]^{\beta}  \ar[d]^{p} &\:E \ar[r]^{\alpha} \ar[d]^{p}  &\:E \ar[r]^{\beta}  \ar[d]^{p} &\:E \ar[r]^{\alpha}  \ar[d]^{p} &\\
\ar[r]^{{\beta}^{\text{tr}}} &\:E\phantom{\rule{0.3ex}{2.4ex}}   \ar[r]^{{\alpha}^{\text{tr}}}&\:E\phantom{\rule{0.3ex}{2.4ex}}    \ar[r]^{{\beta}^{\text{tr}}} &\:E\phantom{\rule{0.3ex}{2.4ex}}    \ar[r]^{{\alpha}^{\text{tr}}} &\:E\phantom{\rule{0.3ex}{2.4ex}}    \ar[r]^{{\beta}^{\text{tr}}} &
  }$$ 
  
 \medskip
 
 It follows that ${\mathrm{Coker}}(\alpha)$ and  ${\mathrm{Coker}}(\beta)$ are dual to each other. If we define
$$\kappa = \begin{pmatrix} 0&-1\\-(v-t)&u+s+\gamma t \end{pmatrix}\qquad \text{and} \qquad \lambda = \begin{pmatrix} 1&0\\ -(u-s) & v-t \end{pmatrix},
$$
then
$$
\kappa \alpha =  \begin{pmatrix} u-s&-(v-t)\\0&0 \end{pmatrix} \qquad \text{and}\qquad \lambda \beta =  \begin{pmatrix} v-t&-(u+s + \gamma t)\\0&0 \end{pmatrix}.
$$

\medskip  

Since $x$ is a monic polynomial in the variable $X$ it follows from the division algorithm that every element of $R$ can be written as a sum $f(v)+ug(v)$, where $f$ and $g$ are uniquely defined polynomials in $A[Y]$. This shows that $v-t$ is not a zero-divisor in $R$, so the maps $\kappa$ and $\lambda$ are injective. Therefore ${\mathrm{Coker}}(\beta)  \approx {\mathrm{Im}}(\alpha)$ is isomorphic to the ideal $J = (u-s,v-t) \subseteq R$, and ${\mathrm{Coker}}(\alpha)  \approx {\mathrm{Im}}(\beta)$ is isomorphic to the ideal $K = (-(v-t), u+s+\gamma t ) \subseteq R$. In the total quotient ring of $R$ multiplication by $v-t$ is an isomorphism, so therefore ${\mathrm{Im}}(\beta)$ is isomorphic to the fractional ideal  $J' =(1, \epsilon)$ where $\epsilon = \frac{u+s+\gamma t}{v-t}$. Tracing all the maps the isomorphism $J' \rightarrow  {\mathrm{Hom}}(J,R) = \dual{J}$ is simply given by multiplication. Note that we have
\[ \epsilon (u-s) = -(\delta v + \delta t + \gamma u)    \quad \text{and} \quad  \epsilon (v-t) = u + s + \gamma t.
\]
Let $J''$ be the fractional ideal consisting of all elements of the total quotient ring multiplying $J$ into $R$. Certainly $J' \subseteq J''$. If $q \in J''$, we can find elements $r,s \in R$ such that $r+ s \epsilon -q$ induces the zero homomorphism on $J$. Then $(r + s \epsilon -q)(v-t) = 0$ and so $q = r  + s \epsilon$ and $J' = J''$.
Again since $v-t$ is a not a zero-divisor, the submodule of  $\dual{J}$ generated by the inclusion in ${\mathrm{Hom}}(J,R)$ is isomorphic to $R$ and the quotient $\dual{J}/R$ is generated by a single element $[\epsilon]$. The map $R \rightarrow \dual{J}/R$ sending the identity to $[\epsilon]$ factors through $A \approx R/(u-s, v-t)$. We claim that the induced map $A \rightarrow \dual{J}/R$ is an isomorphism. To see this, suppose that $r \in A$ and that $r \epsilon \in R$. Then $r(u+s+\gamma t) = (v-t)(f(v) +ug(v))$ where $f$ and  $g$ are uniquely determined. Comparing $v$ terms, we must have $r=f=g=0$. Using \cite{Knudsen:83}, Appendix, Theorem 2, or \cite{Ile:11} (to appear), we summarize these results.
\begin{proposition} \label{prop:key} \quad Given $A$, $\gamma, \delta, s, t, u, v, \alpha, \beta, \epsilon$,  $R$ and $E $ as above, then
\begin{itemize} 
\item[$\bullet$] The ideal $J = (u-s,v-t)$ is stably reflexive. 
\item[$\bullet$] The dual $\dual{J}  = {\mathrm{Hom}}(J,R)$ is isomorphic to the fractional ideal $(J : R)$ and is generated the elements $1$  and $\epsilon = \frac{u+s+\gamma t}{v-t}$.
\item[$\bullet$] The map $A \rightarrow  \dual{J}/R$ defined by multiplication with the class of $\epsilon$ is an isomorphism.
\item[$\bullet$] The sequence \qquad $ \markedrightarrow{\beta} E  \markedrightarrow{\alpha} E  \markedrightarrow{\beta} E  \markedrightarrow{\alpha} E  \markedrightarrow{\beta} E  \markedrightarrow{\alpha} $ \qquad is exact and there are isomorphisms $J \approx {\mathrm{Coker}}(\beta)$ and 
$\dual{J} \approx {\mathrm{Coker}}(\alpha)$.  \hfill  \proofsquare 

\end{itemize} 
\end{proposition}  
\begin{corollary}\label{cor:key}
Let $A$, $\gamma, \delta, s, t, u, v, \alpha, \beta, \epsilon$,  $R$ and $E $ as above, and assume that $A$ is a local ring and that  $s,t \in {\bm{m}}_{A}$ the maximal ideal of $A$,  ${\bm{m}}_{R} = {\bm{m}}_{A}R + (u,v)R$ is a maximal ideal of $R$, and we let ${\widehat{A}}$ and ${\widehat{R}}$ denote the completions of $A$ and $R$. The maps $\alpha$ and $\beta$ extend to maps ${\widehat{\alpha}}$ and ${\widehat{\beta}}$ of ${\widehat{E}}$, and
\begin{itemize} 
\item[$\bullet$] The ideal ${\widehat{J}} = (u-s,v-t){\widehat{R}}$ is stably reflexive. 
\item[$\bullet$] The dual ${\duall{\widehat{J} }} = {\mathrm{Hom}}({\widehat{J}},{\widehat{R}}) \approx {\mathrm{Hom}}(J,R) \otimes_{R}{\widehat{R}}$ is isomorphic to the fractional ideal $(\widehat{J} : \widehat{R})$ and is generated the elements $1$  and $\epsilon = \frac{u+s+\gamma t}{v-t}$.
\item[$\bullet$] The map ${\widehat{A}} \rightarrow {\duall{\widehat{J}}}/{\widehat{R}}$ defined by multiplication with the class of $\epsilon$ is an isomorphism.
\item[$\bullet$] The sequence \qquad $ \markedrightarrow{\widehat{\beta}} {\widehat{E}}  \markedrightarrow{\widehat{\alpha}} {\widehat{E}} \markedrightarrow{\widehat{\beta}} {\widehat{E}} \markedrightarrow{\widehat{\alpha}}  {\widehat{E}}  \markedrightarrow{\widehat{\beta}} {\widehat{E}}  \markedrightarrow{\widehat{\alpha}} 
$ \qquad is exact and there are isomorphisms ${\widehat{J}} \approx {\mathrm{Coker}}({\widehat{\beta}})$ and 
${\duall{\widehat{J}}} \approx {\mathrm{Coker}}({\widehat{\alpha}})$. 
\end{itemize} 
\end{corollary}
\begin{proof} This holds because ${\widehat{R}}$ is flat over $R$.
\end{proof}

\section{Proof of the Main Lemma}
\begin{proof} The statement is of local nature, so let $s \in S$, $x = {\Delta}(s) \in C$. We need only consider the case where $x$ is a double point of $C_{s}= C \times _{S} {\mathrm{Spec}}(k)$, where $k = {\mathcal{O}}_{S,s}/{\bm{m}}_{s}$. We define ${\mathcal{O}}_{C,x}= R$, ${\mathcal{O}}_{S,s}= A$. The section $\Delta$ defines a homomorphism $R \rightarrow A$ which we by abuse of language also call $\Delta$. We have a non-degenerate quadratic form $q(X,Y) = X^{2} + \gamma XY + \delta Y^{2}$ such that ${\widehat{R}} \otimes_{\widehat{A}} k \approx R_{0} = k[[X,Y]]/(q)$. Lifting $\gamma$ and $\delta$ to $\widehat{A}$, it follows from Proposition  \ref{prop:versal} that there exists a map $R_{pd} \rightarrow {{\widehat{\mathcal{O}}_{C,x}}} = {\widehat{R}}$ such that the diagram below is a morphism in $\myscripthat{F}$.
\begin{center}
$
\SelectTips{cm}{} 
\xymatrix{  R_{pd}  \ar[r] \ar@/_3ex /[d]_{{\Delta}_{pd}} &  {\widehat{R}} \ar[r] \ar@/^3ex /[d]^{\widehat{\Delta}}& R_{0} \ar@/^3ex /[d]^{{\Delta}_{0}}  \\ A_{pd} \ar[r] \ar[u]& {\widehat{A}} \ar[r] \ar[u]& k\ar[u]}
$
\end{center}
If $I$ the kernel of $\Delta : R \rightarrow A$, then $\widehat{I}$ is the kernel of $\widehat{\Delta} : \widehat{R} \rightarrow \widehat{A}$. By Corollary \ref{cor:key}, $\widehat{A}$,  $\widehat{R}$ and $\widehat{I}$ satisfy the properties of the Main Lemma, and by  \cite{Knudsen:83}, Appendix, Proposition 6, or \cite{Ile:11}, so do $R$, $A$ and $I$.
\end{proof}

\section{A close look at the fiber}
\begin{theorem}\label{thm:stab}
For every object  $(\sigma : C \rightarrow S, s_{1},\ldots,s_{n},\Delta)$ of ${\overline{\mathcal{C}}}_{g,n}$, $C^{s} = {\textbf{\emph{Proj}}}_{C}({\mathrm{Sym}}\hspace*{1pt}{\mathscr{K}})$ with the sections defined by the surjections $s_{i}^{*}{\mathscr{K}} \rightarrow s_{i}^{*}{\mathscr{L}}_{s}$ and ${\Delta}^{*}{\mathscr{K}} \rightarrow {\Delta}^{*}{\mathscr{L}}_{\Delta}$, is an object of ${\overline{\mathcal{M}}}_{g,n{+}1}$.
\end{theorem}

\begin{proof}\quad Since the case where $\Delta$ meets one of the other sections is satisfactorily treated in \cite{Knudsen:83}, we will only look at the case where $\Delta$ hits a double point of the fiber. We start by analyzing the case we treated in the Key Example. We have the rings $A$ and $R$ together with elements $\gamma$, $\delta$, $s$, $t$, $u$, $v$ where $u$ and $v$  satisfy exactly the one relation
$$(u + s + \gamma t)(u - s)+(\gamma u + \delta v + \delta t)(v - t)=0.
$$
Since
$$\det \begin{pmatrix}  1&0&1&\gamma \\  0&1&0& \delta \\  -1&0&1&0 \\  0&-1&\gamma& \delta \end{pmatrix} = 4\delta - {\gamma}^{2}
$$
is a unit in $A$, we have $(u-s, v-t, u + s + \gamma t, \gamma u + \delta v + \delta t) = (u,v,s,t)$.

\smallskip

Let $p$ be a prime ideal in $R$ containing $I=(u-s,v-t)$. Then $s$ or $t$ is invertible in $R_{p}$ if and only if $u + s + \gamma t$ or $\gamma u + \delta v + \delta t$ is invertible in $R_{p}$.   Since $(u + s + \gamma t)(u-s) = (\gamma u + \delta v + \delta t)(v-t)$ it follows that if $s$ or $t$ is invertible in $R_{p}$, $IR_{p}$ is an invertible ideal, and hence so is $\dual{(IR_{p})}$. 

\smallskip

Let $X= {\mathrm{Spec}}(R)$, $X^{s}= {\mathrm{Proj}}_{R}({\mathrm{Sym}}(\dual{I}))$, $r:X^{s} \rightarrow X$, $\tau : X \rightarrow S = {\mathrm{Spec}}(A)$. The map $\Delta : R \rightarrow A$ define a section which we also name $\Delta : S \rightarrow X$. Let $\mathscr{J}$ be the sheaf of ideals defined by $I$. We have just seen that $\dual{\mathscr{I}}$ is invertible outside the locus $\Delta (V(s,t))$, and therefore the morphism $r$ is an isomorphism here. We now concentrate on the fibers of $r : X^{s} \rightarrow X$ above the locus $\Delta (V(s,t))$.

\smallskip

We have an exact sequence $R^{2}  \markedrightarrow{\alpha} R^{2} \rightarrow \dual{I} \rightarrow 0$. The last arrow is $(1, -\epsilon)$, and 
\[\alpha = \begin{pmatrix}\gamma u + \delta v + \delta t&u + s + \gamma t \\ -(u - s)&v - t \end{pmatrix}.
\]
Therefore $X^{s} = {\mathrm{Proj}}(R[X_{0},X_{1}]/(f,g))$, where $f= (\gamma u + \delta v + \delta t)X_{0} -(u - s)X_{1}$ and $g= (u + s + \gamma t)X_{0} + (v-t)X_{1}$, and the lifted section ${\Delta}^{s}$ sends any point in the locus $V(s,t)$ to the point with homogeneous coordinates $(0,1)$.  We define $x=X_{0}/X_{1}$ and $y=X_{1}/X_{0}$ and $R_{0}= R[y]/(f_{0},g_{0})$ and $R_{1}= R[x]/(f_{1},g_{1})$, where 
\begin{align*}  &f_{0}  =  (\gamma u + \delta v + \delta t) -(u - s)y  \qquad \text{and}\qquad f_{1} =  (\gamma u + \delta v + \delta t)x -(u - s) \\  
                       &g_{0}=  (u + s + \gamma t) + (v-t)y    \quad \qquad \,\text{and}\qquad  g_{1}=  (u + s + \gamma t)x + (v-t) 
\end{align*}
In the ring $R_{0}$ we have the identity $(v-t)f_{0}+(u-s)g_{0}=q(u,v)-q(s,t)=0$, hence all relations between $u$ and $v$ in $R_{0}$ are generated by $f_{0}$ and $g_{0}$ and similarly the relations in $R_{1}$ are generated by $f_{1}$ and $g_{1}$. We may use $g_{0}$ to eliminate $u$ and $g_{1}$ to eliminate $v$. A little arithmetic shows that we have
\begin{align} 
R_{0} &= A[v,y]/(v(y^{2} - \gamma y + \delta) + s(2y - \gamma) + t (- y^{2} +2\gamma y - {\gamma}^{2} + \delta))  \label{eq0:covering} \\ 
R_{1} &= A[u,x]/(u(\delta x^{2} - \gamma x + 1)  + s(\delta x^{2} - 1) + t \delta (\gamma x^{2} - 2x))  \label{eq1:covering}
\end{align} 
We first check flatness over $A$. Every element of $R_{0}$ can be written uniquely as a polynomial in $v$ and $y$ without monomials of the form $v^{n}y^{m}$ with $n \geq 1$ and $m \geq 2$, so it is free as an $A$-module. Since $x=0$ is not a root of $\delta x^{2} - \gamma x + 1=0$, the locus $V(\delta x^{2} - \gamma x + 1)$ is contained in ${\mathrm{Spec}}(R_{0})$, therefore $X^{s}$ is covered by ${\mathrm{Spec}}(R_{0}) \cup {\mathrm{Spec}}(R_{1}[(\delta x^{2} - \gamma x + 1)^{-1}])$. We have $R_{1}[(\delta x^{2} - \gamma x + 1)^{-1}] = A[x][(\delta x^{2} - \gamma x + 1)^{-1}]$ which is flat over $A$ and hence $X^{s}$ is flat over $S$. Outside the locus of $s=t=0$ in $S$, the fibers have not changed. If $s=t=0$, it follows from (\ref{eq0:covering}) and (\ref{eq1:covering}), that the geometric fibers consist of a projective line meeting two parallel affine lines transversally. The image of the lifted section is contained in ${\mathrm{Spec}}(R_{1})$ where it is defined by $x=0$, and here the morphism $\tau r$ is smooth.

\medskip

Back to the general case, let $p: C^{s} \rightarrow C$ and let $x \in S$ be a point such that $y = \Delta(x)$ is a double point of the fiber of $\sigma$. Let $C_{x }= {\sigma}^{-1}(x)=C\times_{S}{\mathrm{Spec}(k)}$.  We have a diagram of cartesian squares and its completion along the various subschemes.

\begin{center}
$
\SelectTips{cm}{} 
\xymatrix{C^{s} \ar[d]_{p} &(C_{x})^{s}   \ar[l]   \approx   C^{s}  \times_{S}{\mathrm{Spec}}(k) \ar[d]^{p \times 1}\\ 
C  \ar[d]_{\sigma}  & C_{x}   \ar[l]   \approx   C     \times_{S}{\mathrm{Spec}}(k) \ar[d]^{\sigma \times 1}\\ 
S    &{\mathrm{Spec}}(k) \ar[l]} 
$ \qquad \qquad \qquad 
$
\SelectTips{cm}{} 
\xymatrix{C^{s} _{/p^{-1}(y)} \ar[d]_{\widehat{p} } & (C_{x})^{s}_{/p^{-1}(y)}   \ar[l]  \ar[d]^{\widehat{p} \times 1} \\ 
C_{/y}  \ar[d]_{\widehat{\sigma}}  & (C_{x})_{/y}   \ar[l]  \ar[d]^{\widehat{\sigma} \times 1} \\ 
S_{/x}    &{\mathrm{Spec}}(k) \ar[l]  }
$
\end{center}
We write $(\widehat{\mathcal{O}}_{S,x}, \mywidehat{\bm{m}}_{x}) = (A, {\bm{m}}_{A})$. For elements $\gamma$, $\delta$ in $A$, $s$, $t$ in ${\bm{m}}_{A}$ let $q(X,Y) = X^{2}+ \gamma  XY + \delta Y^{2}$, and let $R=A[u,v]$ subject to the only relation $q(u,v)=q(s,t)$ and let ${\bm{m}}_{A}R + (u,v)$ be a maximal ideal in $R$.

According to  Proposition \ref{prop:versal}, we can choose $\gamma$, $\delta$ in $A$, $s$, $t$ in ${\bm{m}}_{A}$ such that $q(X,Y) = X^{2}+ \gamma  XY + \delta Y^{2}$ is non-degenerate and 
\[(\widehat{\mathcal{O}}_{C,y}, \mywidehat{\bm{m}}_{y}) = (\widehat{R}, \mywidehat{\bm{m}}_{R})
\]
Let $X= {\mathrm{Spec}}(R)$, and let $z \in X$ correspond to the maximal ideal ${\bm{m}}_{R}$, it follows that the formal schemes ${\mathrm{Spf}}(A)$, $X_{/z}$ and $X^{s}_{/r^{-1}(z)}$  fit into a commutative diagram
\begin{center}
$
\SelectTips{cm}{} 
\xymatrix{C^{s} _{/p^{-1}(y)} \ar[d]_{\widehat{p}}  \ar[r]^{\approx} &X^{s}_{/r^{-1}(z)}  \ar[d]^{\widehat{r}} \\ 
C_{/y}  \ar[d]_{\widehat{\sigma}}   \ar[r]^{\approx}   & X_{/ z}   \ar[d]^{\widehat{\tau}} \\ 
S_{/x}    \ar[r]^{\approx}    & {\mathrm{Spf}}(A) .}
$
\end{center}
Since $\tau r$ is flat, so is ${\widehat{\tau}}{\widehat{r}}$ and therefore so are ${\widehat{\sigma}}{\widehat{p}}$ and $\sigma p$. We have also seen that the geometric fiber at $z$ is a projective line containing the marked point and meeting the rest of the fiber over $x$ transversally in two other distinct points, which means that the fiber is a stable $n{+}1$-pointed curve.
\end{proof}
\begin{theorem}\label{the:stab-funct} The assignment $C$ to $C^{s}$ in Theorem \ref{thm:stab} is a functor $s : {\overline{\mathcal{C}}}_{g,n} \rightarrow {\overline{\mathcal{M}}}_{g,n{+}1}$.
\end{theorem}
\begin{proof}
This follows from \cite{Knudsen:83}, Appendix, Proposition 4. Note that the proof of this proposition is omitted since it is immediate. 
\end{proof}

\begin{center}{\bf{Acknowledgment}}\end{center}
I thank William Fulton for insisting on more accessible proofs of the results in  \cite{Knudsen:83}, Runar Ile for giving me much support and for pointing out to me the particularly important reference  \cite{Eisenbud:80}, Dan Abramovich for reading an earlier version, and not the least the unknown referee for doing a very thorough job.



%% References


\bibliographystyle{model1b-num-names}
\bibliography{<your-bib-database>}

\begin{thebibliography}{00}

%\bibitem{key}...
%%

%\bibitem{}



\bibitem{ACG:11}
Enrico Arbarello, Maurizio Cornalba, Phillip A. Griffiths,
\newblock {\textit{Geometry of Algebraic Curves,  Volume II}},
\newblock Grundlehren der mathematischen Wissenschaften, vol. 268,  963 p., Springer-Verlag 2011.
\quad {\tiny{{\color{MyDarkBlue}\url{http://www.springer.com/mathematics/algebra/book/978-3-540-42688-2}}}}
\bibitem{Eisenbud:80}
David Eisenbud,
\newblock {\textit{Homological Algebra on a Complete Intersection, with an Application to Group Representations}},
\newblock Transactions of the American Mathematical Society, vol. 260, no. 1, (July 1980), 35-64. MR0570778 (82d:13013)
 \quad {\tiny{{\color{MyDarkBlue}\url{http://www.jstor.org/pss/1999875}}}}
\bibitem{Ile:11}
Runar Ile,
\newblock {\textit{Stably Reflexive Modules and a Lemma of Knudsen}},
\newblock
\quad {\tiny{{\color{MyDarkBlue}\url{http://arxiv.org/abs/1110.3909}}}}
\bibitem{Knudsen:83}
Finn Faye Knudsen,
\newblock {\textit{The Projectivity of the Moduli Space of Stable Curves, II: The Stacks $M_{g,n}$}},
\newblock Mathematica Scandinavica, vol. 52, (1983), no. 2, 161-199, MR0702953 (85d:14038a)
\quad {\tiny{{\color{MyDarkBlue}\url{http://www.mscand.dk/article.php?id=2669}}}}
\bibitem{Schlessinger:68}
Michael Schlessinger,
\newblock {\textit{Functors of Artin Rings}},
\newblock Transactions of the American Mathematical Society, vol. 130, no. 2, (Feb., 1968), 208-222, MR0217093 (36 \#184) 
\quad {\tiny{{\color{MyDarkBlue}\url{http://www.jstor.org/pss/1994967}}}}
 



\end{thebibliography}

%% not want to use model1b-num-names.bst.


\end{document}